\documentclass[twocolumn]{autart}

\usepackage{amsmath,amssymb}
\usepackage{graphicx}
\usepackage{savesym}
\savesymbol{AND}
\usepackage{algorithm}
\usepackage{algorithmic}
\usepackage{flushend}
\usepackage{xcolor}

\usepackage{accents}
\newcommand{\ubar}[1]{\underaccent{\bar}{#1}}

\usepackage{bbm}
\allowdisplaybreaks

\theoremstyle{plain}
\newtheorem{ass}{Assumption}
\newtheorem{lemma}{Lemma}

\makeatletter
\newenvironment{proof}[0]{\par
  \normalfont \topsep6\p@\@plus6\p@\relax
  \trivlist
  \item[\hskip\labelsep
        \itshape
    Proof\@addpunct{.}]\ignorespaces
}{%
  \hfill\ \qed\endtrivlist\@endpefalse
}
\makeatother
\makeatletter
\newenvironment{proofof}[1]{\par
  \normalfont \topsep6\p@\@plus6\p@\relax
  \trivlist
  \item[\hskip\labelsep
        \bfseries
    Proof of #1\@addpunct{.}]\ignorespaces
}{%
  \hfill\ \qed\endtrivlist\@endpefalse
}
\makeatother

\DeclareMathOperator{\conv}{conv}
\let\vert\relax \DeclareMathOperator{\vert}{vert}

\newcommand{\bone}{\mathbbm{1}}
\newcommand{\bN}{\mathbb{N}}
\newcommand{\bR}{\mathbb{R}}
\newcommand{\bZ}{\mathbb{Z}}

\newcommand{\cD}{\mathcal{D}}

\newcommand{\cP}{\mathcal{P}}
\newcommand{\cX}{\mathcal{X}}

\newcommand{\LP}{\mathrm{LP}}
\newcommand{\T}{^\top}
\newcommand{\eps}{\varepsilon}
\newcommand{\norm}[1]{\| #1 \|}

\newcommand{\rhob}{\bar{\rho}}
\newcommand{\rhot}{\tilde{\rho}}
\newcommand{\gammab}{\bar{\gamma}}
\newcommand{\gammat}{\tilde{\gamma}}
\newcommand{\xs}{x^\star}

\newcommand{\xbs}{\bar{x}^\star}
\newcommand{\ls}{\lambda^\star}

\newcommand{\lbs}{\bar{\lambda}^\star}
\newcommand{\sub}{\bar{s}}
\newcommand{\slb}{\ubar{s}}
\newcommand{\deltas}{\rho}

\newcommand{\gammaub}{\bar{\gamma}}
\newcommand{\gammalb}{\ubar{\gamma}}

\newcommand{\PbLP}{\overline{\cP}_{\LP}}
\newcommand{\Db}{\overline{\cD}}

\newcommand{\Jsp}{J^\star_{\text{\ref{eq:primal_program}}}}
\newcommand{\Jslp}{J^\star_{\cP_{\LP}}}
\newcommand{\Jslpb}{J^\star_{\PbLP}}

\begin{document}

\begin{frontmatter}
	\runtitle{Multi-agent MILPs: finite-time feasibility and performance guarantees}
	\title{A decentralized approach to multi-agent MILPs:\\ finite-time feasibility and performance guarantees} 

	\thanks{Research was supported by the European Commission under the project UnCoVerCPS, grant number 643921.}
	\thanks[footnoteinfo]{Corresponding author A.~Falsone. Tel. +39-02-23994028. Fax +39-02-23993412.}

	\author[PoliMi]{Alessandro~Falsone}\ead{\\ alessandro.falsone@polimi.it},
	\author[Oxford]{Kostas~Margellos}\ead{kostas.margellos@eng.ox.ac.uk},
	\author[PoliMi]{Maria~Prandini}\ead{maria.prandini@polimi.it}

	\address[PoliMi]{Dipartimento di Elettronica, Informazione e Bioingegneria, Politecnico di Milano, Via Ponzio 34/5, 20133 Milano, Italy}
	\address[Oxford]{Department of Engineering Science, University of Oxford, Parks Road, Oxford, OX1 3PJ, United Kingdom}


	\begin{abstract}
		We address the optimal design of a large scale multi-agent system where each agent has discrete and/or continuous decision variables that need to be set so as to optimize the sum of linear local cost functions, in presence of linear local and global constraints. The problem reduces to a Mixed Integer Linear Program (MILP) that is here addressed according to a decentralized iterative scheme based on dual decomposition, where each agent determines its decision vector by solving a smaller MILP involving its local cost function and constraint given some dual variable, whereas a central unit enforces the global coupling constraint by updating the dual variable based on the tentative primal solutions of all agents.
		An appropriate tightening of the coupling constraint through iterations allows to obtain a solution that is feasible for the original MILP. The proposed approach is inspired by a recent method to the MILP approximate solution via dual decomposition and constraint tightening, and presents the advantage of guaranteeing feasibility in finite-time and providing better performance guarantees.
		The two approaches are compared on a numerical example on plug-in electric vehicles optimal charging.
	\end{abstract}
\end{frontmatter}

\section{Introduction}
In this paper we are concerned with the optimal design of a large-scale system composed of multiple agents, each one characterized by its set of design parameters that should be chosen so as to solve a constrained optimization problem where the agents' decisions are coupled by some global constraint.
More specifically, the goal is to minimize the sum of local linear cost functions, subject to local polyhedral constraints and a global linear constraint.
A key feature of our framework is that design parameters can have both continuous and discrete components.

Let $m$ denote the number of agents. Then, the optimal design problem takes the form of the following Mixed Integer Linear Program (MILP):
\begin{align}
	\min_{x_1,\dots,x_m} \qquad &\sum_{i=1}^m c_i\T x_i \tag{$\cP$} \label{eq:primal_program} \\
		\text{subject to:} \quad &\sum_{i=1}^m A_i x_i \leq b \nonumber \\
							&x_i \in X_i, \, \, i=1,\dots,m \nonumber
\end{align}
where, for all $i=1,\dots,m$, $x_i\in\bR^{n_i}$ is the decision vector of agent $i$, $c_i\T x_i$ its local cost, and $X_i = \{ x_i\in\bR^{n_{c,i}}\times\bZ^{n_{d,i}}: D_i x_i \leq d_i \}$ its local constraint set defined by a matrix $D_i$ and a vector $d_i$ of appropriate dimensions, $n_{c,i}$ being the number of continuous decision variables and $n_{d,i}$ the number of discrete ones, with $n_{c,i}+n_{d,i}=n_i$.
The coupling constraint $\sum_{i=1}^m A_i x_i \leq b$ is defined by matrices $A_i\in\bR^{p}\times \bR^{n_i}$, $i=1,\dots,m$, and a $p$-dimensional vector $b\in\bR^{p}$.

Despite the advances in numerical methods for integer optimization, when the number of agents is large, the presence of {discrete} decision variables makes the optimization problem hard to solve, and calls for some decomposition into lower scale MILPs, as suggested in \cite{vujanic2016decomposition}.

A common practice to handle problems of the form of \ref{eq:primal_program} consists in first dualizing the coupling constraint introducing a vector $\lambda\in\bR^{p}$ of $p$ Lagrange multipliers and solving the dual program
\begin{equation}
	\max_{\lambda\geq 0} \; -\lambda\T b + \sum_{i=1}^m \min_{x_i\in X_i} (c_i\T + \lambda\T A_i)x_i, \tag{$\cD$} \label{eq:dual_program}
\end{equation}
to obtain $\ls$, and then constructing a primal solution $x(\ls) = [x_1(\ls)\T\;\cdots\;x_m(\ls)\T]\T$ by solving $m$ MILPs given by:
\begin{equation} \label{eq:primal_from_lambda}
	x_i(\lambda) \in \arg\min_{x_i\in \vert(X_i)} (c_i\T + \lambda\!\T A_i)x_i,
\end{equation}
where the search within the closed constraint polyhedral set $X_i$ can be confined to its set of vertices $\vert(X_i)$ since the cost function is linear.

Unfortunately, while this procedure guarantees $x(\ls)$ to satisfy the local constraints since $x_i(\ls)\in X_i$ for all $i=1,\dots,m$, it does not guarantees the satisfaction of the coupling constraint.

A way to enforce the satisfaction of the coupling constraint is to solve \ref{eq:dual_program} via a subgradient method, and then use a recovery procedure for the primal variables, \cite{shor1985minimization}. Albeit this method is very useful in applications since it allows for a distributed implementation, see e.g. \cite{falsone2017automatica,simonetto2016primal}, it provides a feasible solution only when there are no discrete decision variables.
As a matter of fact, if we let $\conv(X_i)$ denote the convex hull of all points inside $X_i$, then, the primal solution recovered using \cite{shor1985minimization} is the optimal solution $\xs_{\LP}$ of the following Linear Program (LP):
\begin{align}
	\min_{x_1,\dots,x_m} \qquad &\sum_{i=1}^m c_i\T x_i \tag{$\cP_{\LP}$} \label{eq:LP_primal_program} \\
		\text{subject to:} \quad &\sum_{i=1}^m A_i x_i \leq b \nonumber \\
							&x_i \in \conv(X_i), \, \, i=1,\dots,m. \nonumber
\end{align}
The dual of the convexified problem \ref{eq:LP_primal_program} coincides with the one of \ref{eq:primal_program} and is given by \ref{eq:dual_program} (see \cite{geoffrion1974lagrangean} for a proof). Clearly $\xs_{\LP} \in \conv(X_1)\times \dots \times \conv(X_m)$ does not necessarily imply that $\xs_{\LP} \in X_1\times X_2 \times \dots \times X_m$. Therefore the solution $\xs_{\LP}$ recovered using \cite{shor1985minimization} satisfies the coupling constraint but not necessarily the local constraints.

For these reasons recovery procedures for MILPs are usually composed of two steps: a tentative solution that is not feasible for either the joint constraint or the local ones is first obtained through duality, and then a heuristic is applied to recover feasibility starting from this tentative solution, see, e.g., \cite{bertsekas1983optimal,redondo1999short}.

Problems in the form of \ref{eq:primal_program} arise in different contexts like power plants generation scheduling \cite{yamin2004review} where the agents are the generation units with their on/off state modeled with binary variables and the joint constraint consists in energy balance equations, or buildings energy management \cite{ioli2015iterative}, where the cost function is a cost related to power consumption and constraints are related to capacity, comfort, and actuation limits of each building. Other problems that fits the structure of \ref{eq:primal_program} are supply chain management \cite{dawande2006effective}, portfolio optimization for small investors \cite{baumann2013portfolio}, and plug-in electric vehicles \cite{vujanic2016decomposition}. In all these cases it is of major interest to guarantee that the derived (primal) solution is implementable in practice, which means that it must be feasible for \ref{eq:primal_program}.

Interestingly, a large class of dynamical systems involving both continuous and logic components can be modeled as a Mixed Logical Dynamical (MLD) system, \cite{Bemporad:1999}, which are described by linear equations and inequalities involving both discrete and continuous inputs and state variables. Model predictive control problems for MLD systems involving the optimization of a linear finite-horizon cost function then also fit the MILP description in \ref{eq:primal_program}.

\subsection*{Background}
Problems in the form of \ref{eq:primal_program} have been investigated in \cite{aubin1976estimates}, where the authors studied the behavior of the duality gap (i.e. the difference between the optimal value of \ref{eq:primal_program} and \ref{eq:dual_program}) showing that it decreases relatively to the optimal value of \ref{eq:primal_program} as the number of agents grows. The same behavior has been observed in \cite{bertsekas1983optimal}. In the recent paper \cite{vujanic2016decomposition}, the authors explored the connection between the solutions $\xs_{\LP}$ to the linear program \ref{eq:LP_primal_program} and $x(\ls)$ recovered via \eqref{eq:primal_from_lambda} from the solution $\ls$ to the dual program \ref{eq:dual_program}. They proposed a method to recover a primal solution which is feasible for \ref{eq:primal_program} by using the dual optimal solution of a modified primal problem, obtained by tightening the coupling constraint by an appropriate amount.

Let $\rho\in \bR^p$ with $\rho\geq0$ and consider the following pair of primal-dual problems:
\begin{align}
	\min_{x_1,\dots,x_m} \qquad &\sum_{i=1}^m c_i\T x_i \tag{${\cP}_{\LP,\rho}$} \label{eq:LP_primal_program_rho} \\
		\text{subject to:} \quad &\sum_{i=1}^m A_i x_i \leq b-\rho \nonumber \\
							&x_i \in \conv(X_i), \, \, i=1,\dots,m \nonumber
\end{align}
and
\begin{equation}
	\max_{\lambda\geq 0} \; -\lambda\T (b-\rho) + \sum_{i=1}^m \min_{x_i\in X_i} (c_i\T + \lambda\T A_i)x_i. \tag{${\cD}_\rho$} \label{eq:dual_program_rho}
\end{equation}
\ref{eq:LP_primal_program_rho} constitutes a tightened version of \ref{eq:LP_primal_program}, whereas \ref{eq:dual_program_rho} is the corresponding dual. For all $j=1,\dots,p$, let $\rhot\in \bR^p$ be defined as follows:
\begin{equation} \label{eq:paul_rho}
	[\rhot]_j = p\max_{i\in\{ 1,\dots, m\}}\ \Big\{ \max_{x_i\in X_i} [A_i]_j x_i - \min_{x_i\in X_i} [A_i]_j x_i \Big\},
\end{equation}
where $[A_i]_j$ denotes the $j$-th row of $A_i$ and $[\rhot]_j$ the $j$-th entry of $\rhot$.

\begin{ass}[Uniqueness] \label{ass:uniqueness_paul}
	Problems \ref{eq:LP_primal_program_rho} and \ref{eq:dual_program_rho} with $\rho$ set equal to $\rhot$ defined in \eqref{eq:paul_rho} have unique solutions $\xs_{\LP,\rhot}$ and $\ls_{\rhot}$, respectively.
\end{ass}
From \cite{vujanic2016decomposition}, we have the following result:
\begin{prop}[Theorem~3.1 in \cite{vujanic2016decomposition}] \label{prop:paul_result}
	Let $\ls_{\rhot}$ be the solution to \ref{eq:dual_program_rho} with $\rho=\rhot$ given in \eqref{eq:paul_rho}. Under Assumption~\ref{ass:uniqueness_paul}, we have that any $x(\ls_{\rhot})$ satisfying \eqref{eq:primal_from_lambda}, is feasible for \ref{eq:primal_program}.
\end{prop}

Let us define
\begin{equation} \label{eq:paul_gamma}
	\gammat = p \max_{i\in\{ 1,\dots, m\}} \left\{ \max_{x_i\in X_i} c_i\T x_i - \min_{x_i\in X_i} c_i\T x_i \right\}.
\end{equation}
Consider the following assumption:
\begin{ass}[Slater] \label{ass:slater_paul}
	There exist a scalar $\zeta > 0$ and $\hat{x}_i\in\conv(X_i)$ for all $i=1,\dots,m$, such that $\sum_{i=1}^m A_i \hat{x}_i \leq b - \rhot - m\zeta\bone$, where $\bone\in\bR^p$ is a vector whose elements are equal to one.
\end{ass}
Then, the sub-optimality level of the approximate solution $x(\ls_{\rhot})$ to \ref{eq:primal_program} can be quantified as follows:
\begin{prop}[Theorem~3.3 in \cite{vujanic2016decomposition}] \label{prop:paul_perf_result}
	Let $\ls_{\rhot}$ be the solution to \ref{eq:dual_program_rho} with $\rho=\rhot$. Under Assumptions~\ref{ass:uniqueness_paul} and~\ref{ass:slater_paul}, we have that $x(\ls_{\rhot})$ derived from \eqref{eq:primal_from_lambda} with $\lambda=\ls_{\rhot}$ satisfies
	\begin{equation} \label{eq:paul_perf_result}
		\sum_{i=1}^m c_i\T x_i(\ls_{\rhot}) - \Jsp \leq \gammat + \frac{\norm{\rhot}_\infty}{p\zeta}\gammat,
	\end{equation}
where $\Jsp $ is the optimal cost of \ref{eq:primal_program}.
\end{prop}

Note that both Proposition~\ref{prop:paul_result} on feasibility and Proposition~~\ref{prop:paul_perf_result} on optimality require the knowledge of the dual solution $\ls_{\rhot}$. This may pose some issues if $\ls_{\rhot}$ cannot be computed centrally, which is the case, e.g., when the agents are not willing to share with some central entity their private information coded in their local cost and constraint set. In those cases, the value of $\ls_{\rhot}$ can only be achieved asymptotically using a decentralized/distributed scheme to solve \ref{eq:dual_program_rho} with $\rho=\rhot$.

\subsection*{Contribution of this paper}

In this paper we propose a decentralized iterative procedure which provides in a finite number of iterations a solution that is feasible for the optimal design problem \ref{eq:primal_program}, thus overcoming the issues regarding the finite-time computability of a decentralized solution in \cite{vujanic2016decomposition}.
Furthermore, the performance guarantees quantifying the sub-optimality level of our solution with respect to the optimal one of \ref{eq:primal_program} are less conservative than those derived in \cite{vujanic2016decomposition}.

As in the inspiring work in \cite{vujanic2016decomposition}, we still exploit some tightening of the coupling constraint to enforce feasibility. However, the amount of tightening is decided through the iterations, based on the explored candidate solutions $x_i\in X_i$, $i=1,\dots,m$, and not using the overly conservative worst-case tightening \eqref{eq:paul_rho} as in \cite{vujanic2016decomposition} where for all $i=1,\dots,m$, the $\max$ and $\min$ of $[A_i]_j x_i$ are computed letting $x_i$ vary over the whole set $X_i$.
The amount of tightening plays a crucial role in the applicability of Proposition~\ref{prop:paul_result}. In fact, a too large value of $\rhot$ may prevent \ref{eq:LP_primal_program_rho} to be feasible when $\rho$ is set equal to $\rhot$, thus violating Assumption~\ref{ass:uniqueness_paul}. A less conservative way to select an appropriate amount of tightening can extend the applicability of the approach to a larger class of problems.
According to a similar reasoning, we are able to improve the bound on the performance degradation of our solution with respect to the optimal one of \ref{eq:primal_program} by taking a less conservative value for the quantity $\gammat$ in \eqref{eq:paul_gamma} that is used in the performance bound \eqref{eq:paul_perf_result}.

Notably, the proposed decentralized scheme allows agents to preserve the privacy on their local information, since they do not have to send to the central unit either their cost coefficients or their local constraints.

\section{Proposed approach}
We next introduce Algorithm~\ref{algo:Alg1} for the decentralized computation in a finite number of iterations of an approximate solution to \ref{eq:primal_program} that is feasible and improves over the solution in \cite{vujanic2016decomposition} both in terms of amount of tightening and performance guarantees.

\begin{algorithm}[t]
	\begin{algorithmic}[1]
		\STATE $\lambda(0) = 0$
		\STATE $\sub_i(0) = -\infty$, $i=1,\dots,m$
		\STATE $\slb_i(0) = +\infty$, $i=1,\dots,m$
		\STATE $k = 0$
		\REPEAT
			\FOR{$i=1$ \TO $m$}
				\STATE $x_i(k+1) \gets \displaystyle\arg\min_{x_i\in \vert(X_i)} (c_i\T + \lambda(k)\T A_i)x_i$ \label{step:primal_update}
			\ENDFOR
			\STATE $\sub_i(k+1) = \max \{ \sub_i(k), A_i x_i(k+1) \}$, $i=1,\dots,m$ \label{step:upper_bounds_update}
			\STATE $\slb_i(k+1) = \min \{ \slb_i(k), A_i x_i(k+1) \}$, $i=1,\dots,m$ \label{step:lower_bounds_update}
			\STATE $\deltas_i(k+1)=\sub_i(k+1) - \slb_i(k+1)$, $i=1,\dots,m$ \label{step:delta_bounds_update}
			\STATE $\rho(k+1) = \displaystyle p \max \{ \deltas_1(k+1), \dots, \deltas_m(k+1) \}$ \label{step:rho_update}
			\STATE $\lambda(k+1)$\\
				$ = \Big[ \lambda(k) + \alpha(k) \big(\displaystyle \sum_{i=1}^m A_i x_i(k+1) - b + \rho(k+1) \big) \Big]_+$ \label{step:dual_update}
			\STATE $k \gets k+1$
		\UNTIL{some stopping criterion is met.}
	\end{algorithmic}
	\caption{Decentralized MILP}
	\label{algo:Alg1}
\end{algorithm}

Algorithm~\ref{algo:Alg1} is a variant of the dual subgradient algorithm. As the standard dual subgradient method, it includes two main steps: step~\ref{step:primal_update} in which a subgradient of the dual objective function is computed by fixing the dual variables and minimizing the Lagrangian with respect to the primal variables, and step~\ref{step:dual_update} which involves a dual update step with step size equal to $\alpha(k)$, and a projection onto the non-negative orthant (in Algorithm~\ref{algo:Alg1} $[\,\cdot\,]_+$ denotes the projection operator onto the $p$-dimensional non-negative orthant $\bR^p_+$).
The operators $\max$ and $\min$ appearing in steps \ref{step:upper_bounds_update}, \ref{step:lower_bounds_update}, and \ref{step:delta_bounds_update} of Algorithm~\ref{algo:Alg1} with arguments in $\bR^p$ are meant to be applied component-wise.
The sequence $\{\alpha(k)\}_{k\geq 0}$ is chosen so as to satisfy $\lim_{k\to\infty} \alpha(k) = 0$ and $\sum_{k=0}^\infty \alpha(k) = \infty$, as requested in the standard dual subgradient method to achieve asymptotic convergence.
Furthermore, in order to guarantee that the solution to step~\ref{step:primal_update} in Algorithm~\ref{algo:Alg1} is well-defined, we impose the following assumption on \ref{eq:primal_program}:
\begin{ass}[Boundedness] \label{ass:boundedness}
	The polyhedral sets $X_i$, $i=1,\dots,m$, in problem \ref{eq:primal_program} are bounded.
\end{ass}

If $\arg\min_{x_i\in \vert(X_i)} (c_i\T + \lambda(k)\T A_i)x_i$ in step~\ref{step:primal_update} is a set of cardinality larger than 1, then, a deterministic tie-break rule is applied to choose a value for $x_i(k+1)$.

Algorithm~\ref{algo:Alg1} is conceived to be implemented in a decentralized scheme where, at each iteration $k$, every agent $i$ updates its local tentative solution $x_i(k+1)$ and communicates $A_ix_i(k+1)$ to some central unit that is in charge of the update of the dual variable. The tentative value $\lambda(k+1)$ for the dual variable is then broadcast to all the agents. Note that the agents do not need to communicate to the central unit their private information regarding their local constraint set and cost but only their tentative solution $x_i(k)$.

The tentative primal solutions $x_i(k+1)$, $i=1,\dots,m$, computed at step~\ref{step:primal_update} are used in Algorithm~\ref{algo:Alg1} by the central unit to determine the amount of tightening $\rho(k+1)$ entering step~\ref{step:dual_update}. The value of $\rho(k+1)$ is progressively refined through iterations based only on those values of $x_i\in X_i$, $i=1, \dots,m$, that are actually considered as candidate primal solutions, and not based on the whole sets $X_i$, $i=1, \dots,m$. This reduces conservativeness in the amount of tightening and also in the performance bound of the feasible, yet suboptimal, primal solution.

Algorithm~\ref{algo:Alg1} terminates after a given stopping criteria is met at the level of the central unit, e.g., if for a given number of subsequent iterations $x(k) = [x_1(k)\T\;\cdots\;x_m(k)\T]\T$ satisfies the coupling constraint. As shown in the numerical study in Section~\ref{sec:exaple}, variants of Algorithm~\ref{algo:Alg1} can be conceived to get an improved solution in the same number of iterations of Algorithm~\ref{algo:Alg1}. The agents should however share with the central entity additional information on their local cost, thus partly compromising privacy preservation.

As for the initialization of Algorithm~\ref{algo:Alg1}, $\lambda(0)$ is set equal to $0$ so that at iteration $k=0$ each agent $i$ computes its locally optimal solution
\begin{equation*}
	x_i(1) \gets \arg\min_{x_i\in \vert(X_i)} c_i\T x_i.
\end{equation*}
Since $\rho(1) = 0$, if the local solutions $x_i(1)$, $i=1, \dots,m$, satisfy the coupling constraint (and they hence are optimal for the original problem \ref{eq:primal_program}), then, Algorithm~\ref{algo:Alg1} will terminate since $\lambda$ will remain 0, and the agents will stick to their locally optimal solutions.


Before stating the feasibility and performance guarantees of the solution computed by Algorithm~\ref{algo:Alg1}, we need to introduce some further quantities and assumptions.

Let us define for any $k\ge 1$
\begin{align}
	\gamma(k) = p \max_{i\in\{1,\dots,m\}} \Big\{ \max_{r\le k} c_i\T x_i(r) - \min_{r\le k} c_i\T x_i(r) \Big\},\label{eq:gamma_update}
\end{align}
where $\{x_i(r)\}_{r\ge 1}$, $i=1,\dots,m$, are the tentative primal solutions computed at step~\ref{step:primal_update}.

Due to Assumption~\ref{ass:boundedness}, for any $i=1,\dots,m$, $\conv(X_i)$ is a bounded polyhedron. If it is also non-empty, then $\vert(X_i)$ is a non-empty finite set (see Corollaries~2.1 and 2.2 together with Theorem~2.3 in \cite[Chapter 2]{bertsimas1997introduction}). As a consequence, the sequence $\{\gamma(k)\}_{k\geq 1}$ takes values in a finite set. Since this is a monotonically non-decreasing sequence, it converges in finite-time to some value $\gammab$.
\\
The same reasoning can be applied to show that the sequence $\{\rho(k)\}_{k\geq 1}$, iteratively computed in Algorithm~\ref{algo:Alg1} (see step~\ref{step:rho_update}), and given by
\begin{align*}
	[\rho(k)]_j = p \max_{i\in\{1,\dots,m\}} \Big\{ \max_{r\le k} [A_i]_j x_i(r) - \min_{r\le k} [A_i]_j x_i(r) \Big\},
\end{align*}
for $j=1,\dots,p$, converges in finite-time to some $\rhob$ since it takes values in a finite set and is (component-wise) monotonically non-decreasing. Note that the limiting values $\rhob$ and $\gammab$ for $\{\rho(k)\}_{k\geq 1}$ and $\{\gamma(k)\}_{k\geq 1}$ satisfy $\rhob\leq\rhot$ and $\gammab\leq\gammat$ where $\rhot$ and $\gammat$ are defined in \eqref{eq:paul_rho} and \eqref{eq:paul_gamma}.

Similarly to \cite{vujanic2016decomposition}, define $\PbLP$ and $\Db$ as the primal-dual pair of optimization problems that are given by setting $\rho$ equal to $\rhob$ in \ref{eq:LP_primal_program_rho} and \ref{eq:dual_program_rho}.

\begin{ass}[Uniqueness] \label{ass:uniqueness}
	Problems $\PbLP$ and $\Db$ have unique solutions $\xbs_{\LP}$ and $\lbs$.
\end{ass}
\begin{ass}[Slater] \label{ass:slater}
	There exists a scalar $\zeta > 0$ and $\hat{x}_i\in\conv(X_i)$ for all $i=1,\dots,m$, such that $\sum_{i=1}^m A_i \hat{x}_i \leq b - \rhob - m\zeta\bone$.
\end{ass}
Note that, since $\rhob\leq\rhot$, if Assumption~\ref{ass:slater_paul} is satisfied, then Assumption~\ref{ass:slater} is automatically satisfied.

We are now in a position to state the two main results of the paper.

\begin{thm}[Finite-time feasibility] \label{thm:primal_feasibility}
	Under Assumptions~\ref{ass:boundedness} and~\ref{ass:uniqueness}, there exists a finite iteration index $K$ such that, for all $k\geq K$, $x(k)=[x_1(k)\T\;\cdots\;x_m(k)\T]\T$, where $x_i(k)$, $i=1,\dots,m$, are computed by Algorithm~\ref{algo:Alg1}, is a feasible solution for problem \ref{eq:primal_program}, i.e., $\sum_{i=1}^m A_i x_i(k) \leq b$, $k\ge K$ and $x_i(k)\in X_i$, $i=1,\dots,m$.
\end{thm}

\begin{thm}[Performance guarantees] \label{thm:primal_performance}
	Under Assumptions~\ref{ass:boundedness}-\ref{ass:slater}, there exists a finite iteration index $K$ such that, for all $k\geq K$, $x(k)=[x_1(k)\T\;\cdots\;x_m(k)\T]\T$, where $x_i(k)$, $i=1,\dots,m$, are computed by Algorithm~\ref{algo:Alg1}, is a feasible solution for problem \ref{eq:primal_program} that satisfies the following performance bound:
	\begin{equation} \label{eq:perf_result}
		\sum_{i=1}^m c_i\T x_i(k) - \Jsp \leq \gammab + \frac{\norm{\rhob}_\infty}{p\zeta}\gammat.
	\end{equation}
\end{thm}
By a direct comparison of \eqref{eq:paul_perf_result} and \eqref{eq:perf_result} we can see that the bound in \eqref{eq:perf_result} is no worse than \eqref{eq:paul_perf_result} due to the fact that $\rhob\leq\rhot$ and $\gammab\leq\gammat$.

\section{Proof of the main results}

\subsection{Preliminary results}

\begin{prop}[Dual asymptotic convergence] \label{prop:dual_optimality}
	Under Assumptions~\ref{ass:boundedness} and~\ref{ass:uniqueness}, the Lagrange multiplier sequence $\{\lambda(k)\}_{k\geq 0}$ generated by Algorithm~\ref{algo:Alg1} converges to an optimal solution of $\Db$.
	\begin{proof}
		As discussed after equation \eqref{eq:gamma_update}, there exists a $K\in\bN$ such that for all $k\geq K$ we have that the tightening coefficient $\rho(k)$ computed in Algorithm~\ref{algo:Alg1} becomes constant and equal to $\rhob$. Therefore, for any $k\geq K$, Algorithm~\ref{algo:Alg1} reduces to the following two steps
		\begin{align}
			&x_i(k+1) \in \arg\min_{x_i\in \vert(X_i)} (c_i\T + \lambda(k)\T A_i)x_i \label{eq:primal_update_afterK} \\
			&\lambda(k+1) = {\left[ \lambda(k) + \alpha(k) \left(\sum_{i=1}^m A_i x_i(k+1) - b + \rhob \right) \right]}_+ \label{eq:dual_update_afterK}
		\end{align}
		which constitute a gradient ascent iteration for $\Db$. According to \cite{bertsekas1999nonlinear}, the sequence $\{\lambda(k)\}_{k\geq 0}$ generated by the iterative procedure \eqref{eq:primal_update_afterK}-\eqref{eq:dual_update_afterK} is guaranteed to converge to the (unique under Assumption \ref{ass:uniqueness}) optimal solution of $\Db$.
	\end{proof}
\end{prop}

\begin{lemma}[Robustness against cost perturbation] \label{lemma:LP_sensitivity}
	Let $P$ be a non-empty bounded polyhedron. Consider the linear program $\min_{x\in P} (c\T+\delta\T)x$, where $\delta$ is a perturbation in the cost coefficients. Define the set of optimal solutions as $\cX(\delta)$.
	There always exists an $\eps > 0$ such that for all $\delta$ satisfying $\norm{\delta} < \eps$, we have $\cX(\delta) \subseteq \cX(0)$.
	\begin{proof}
		Let $u(\delta) = \min_{x\in P} (c\T+\delta\T)x$. Since $P$ is a bounded polyhedron, the minimum is always attained and $u(\delta)$ is finite for any value of $\delta$. The set $\cX(\delta)$ can be defined as
		\begin{equation} \label{eq:lemma_opt_set}
			\cX(\delta) = \{x\in P: (c\T+\delta\T)x \le u(\delta)\},
		\end{equation}
		which is a non-empty polyhedron. As such, it can be described as the convex hull of its vertices (see Theorem~2.9 in \cite[Chapter 2]{bertsimas1997introduction}), which are also vertices of $P$ (Theorem~2.7 in \cite[Chapter 2]{bertsimas1997introduction}).
		
		Let $V = \vert(P)$ and $V_\delta = \vert(\cX(\delta))\subseteq V$. Consider $\delta = 0$.
		
		If $V_0=V$, then, given the fact that, for any $\delta$, $\cX(\delta)$ is the convex hull of $V_\delta$ and $V_\delta \subseteq V=V_0$, we have trivially that $\cX(\delta) \subseteq \cX(0)$, for any $\delta$.
		
		Suppose now that $V_0\subset V$. For any choice of $\xs\in V_0$ and $x\in V\setminus V_0$, we have that $c\T \xs < c\T x$, or equivalently $c\T (\xs-x) < 0$. Pick
		\begin{equation} \label{eq:lemma_eps_def}
			\eps = \min_{\substack{\xs\in V_0 \\ x \in V\setminus V_0}} -\frac{c\T (\xs-x)}{\norm{\xs-x}}
		\end{equation}
		and let $(\bar{x}^\star,\bar{x})$ be the corresponding minimizer. By construction, \eqref{eq:lemma_eps_def} is well defined since $\bar{x}^\star$ is different from $\bar{x}$. Since $c\T (\xs-x) < 0$ for any $\xs\in V_0$ and $x\in V\setminus V_0$, we have that $\eps > 0$. Moreover, for any $\xs\in V_0$ and $x\in V\setminus V_0$, if $\delta$ satisfies $\norm{\delta} < \eps$, then
		\begin{align}
			(c\T+\delta\T)&(\xs-x) = c\T(\xs-x) + \delta\T(\xs-x) \nonumber \\
				&\leq c\T(\xs-x) + \norm{\delta} \norm{\xs-x} \nonumber \\
				&< c\T(\xs-x) + \eps \norm{\xs-x} \nonumber \\
				&\leq c\T(\xs-x) + \left( -\frac{c\T(\xs-x)}{\norm{\xs-x}} \right) \norm{\xs-x} \nonumber \\
				&= c\T(\xs-x) -c\T(\xs-x) = 0, \label{eq:lemma_opt_ineq}
		\end{align}
		where the first inequality is given by the fact that $u\T v \leq |u\T v|$ together with the Cauchy--Schwarz inequality $|u\T v| \leq \norm{u}\norm{v}$, the second inequality is due to $\delta$ satisfying $\norm{\delta} < \eps$, and the third inequality is given by the definition of $\eps$ in \eqref{eq:lemma_eps_def}. \\
		By \eqref{eq:lemma_opt_set} and the definition of $u(\delta)$, for any point $x_\delta$ in the set $V_\delta$, we have that $(c\T+\delta\T)x_\delta \leq (c\T+\delta\T)x$, for all $x\in V$, and therefore $(c\T+\delta\T)x_\delta \leq (c\T+\delta\T)\xs$ for any $\xs\in V_0\subset V$.
		By \eqref{eq:lemma_opt_ineq}, whenever $\norm{\delta}<\eps$, we have that $(c\T+\delta\T)\xs < (c\T+\delta\T)x$ for any choice of $\xs\in V_0$ and $x\in V\setminus V_0$, therefore $(c\T+\delta\T)x_\delta < (c\T+\delta\T)x$ for any $x\in V\setminus V_0$. Since the inequality is strict, we have that $x_\delta \not\in V\setminus V_0$, which implies $x_\delta \in V_0$. Since this holds for any $x_\delta\in V_\delta$, we have that $V_\delta \subseteq V_0$. \\
		Finally, given the fact that, for any $\delta$, $\cX(\delta)$ is the convex hull of $V_\delta$ and $V_\delta \subseteq V_0$, we have $\cX(\delta) \subseteq \cX(0)$, thus concluding the proof.
	\end{proof}
\end{lemma}

Exploiting Lemma~\ref{lemma:LP_sensitivity}, we shall show next that each $\{x_i(k)\}_{k\geq 1}$ sequence, $i=1,\dots,m$, converges in finite-time to some set.

\begin{prop}[Primal finite-time set convergence] \label{prop:finite_time_primal_set}
	Under Assumptions~\ref{ass:boundedness} and~\ref{ass:uniqueness}, there exists a finite $K$ such that for all $i=1,\dots,m$ the tentative primal solution $x_i(k)$ generated by Algorithm~\ref{algo:Alg1} satisfies
	\begin{equation} \label{eq:primal_solution_afterKeps}
		x_i(k) \in \arg\min_{x_i\in \vert(X_i)} (c_i\T + \lbs\!\T A_i)x_i,\quad k\geq K,
	\end{equation}
	where $\lbs$ is the limit value of the Lagrange multiplier sequence $\{\lambda(k)\}_{k\geq 0}$.
	\begin{proof}
		Consider agent $i$, with $i \in \{1,\dots,m\}$. We can characterize the solution $x_i(k)$ in step~\ref{step:primal_update} of Algorithm~\ref{algo:Alg1} by performing the minimization over $\conv(X_i)$ instead of $\vert(X_i)$ since the problem is linear and by enlarging the set $\vert(X_i)$ to $\conv(X_i)$ we still obtain all minimizers that belong to $\vert(X_i)$. Adding and subtracting $\lbs\!\T A_i x_i$ to the cost, we then obtain
		\begin{equation} \label{eq:primal_update_modified}
			x_i(k) \in \arg\min_{x_i\in \conv(X_i)} (c_i\T + \lbs\!\T A_i + (\lambda(k-1)-\lbs)\T A_i)x_i.
		\end{equation}
		
		Set $\delta_i(k-1)\T = (\lambda(k-1)-\lbs)\T A_i$, and let $\cX_i(\delta_i(k-1))$ be the set of minimizers of \eqref{eq:primal_update_modified} as a function of $\delta_i(k-1)$. By Lemma~\ref{lemma:LP_sensitivity}, we know that there exists an $\eps_i > 0$ such that if $\norm{\delta_i(k-1)} < \eps_i$, then $\cX_i(\delta_i(k-1)) \subseteq \cX_i(0)$.
		
		Since, by Proposition~\ref{prop:dual_optimality}, the sequence $\{\lambda(k)\}_{k\geq 0}$ generated by Algorithm~\ref{algo:Alg1} converges to $\lbs$, by definition of limit, we know that there exists a $K_i$ such that $\norm{\delta_i(k-1)} = \norm{(\lambda(k-1)-\lbs)\T A_i} < \eps_i$ for all $k\geq K_i$. Therefore, for every $k\geq K=\max\{K_1, \dots, K_m\}$, we have that $x_i(k) \in \cX_i(0)=\arg\min_{x_i\in \conv(X_i)} (c_i\T + \lbs\!\T A_i)x_i$, $i=1,\dots,m$. This property jointly with the fact that $x_i(k) \in\vert(X_i)$, $i=1,\dots,m$, leads to \eqref{eq:primal_solution_afterKeps}, thus concluding the proof.
	\end{proof}
\end{prop}

\subsection{Proof of Theorems~\ref{thm:primal_feasibility} and~\ref{thm:primal_performance}}

\begin{proofof}{Theorem~\ref{thm:primal_feasibility}}
	Theorem~2.5 of \cite{vujanic2016decomposition} establishes a relation between the solution $\xbs_{\LP}$ of $\PbLP$ and the one recovered in \eqref{eq:primal_from_lambda} from the optimal solution $\lbs$ of the dual optimization problem $\Db$. Specifically, it states that there exists a set of indices $I\subseteq\{1,\dots,m\}$ of cardinality at least $m-p$, such that $[\xbs_{\LP}]^{(i)} = x_i(\lbs)$ for all $i\in I$, where $[\xbs_{\LP}]^{(i)}$ is the subvector of $\xbs_{\LP}$ corresponding to the $i$-th agent. Therefore, following the proof of Theorem~3.1 in \cite{vujanic2016decomposition}, we have that
	\begin{align}
		\sum_{i=1}^m &A_i x_i(\lbs) \nonumber \\
			&= \sum_{i\in I} A_i x_i(\lbs) + \sum_{i\in I^c} A_i x_i(\lbs) \nonumber \\
			&= \sum_{i\in I} A_i [\xbs_{\LP}]^{(i)} + \sum_{i\in I^c} A_i x_i(\lbs) \nonumber \\
			&= \sum_{i=1}^m A_i [\xbs_{\LP}]^{(i)} + \sum_{i\in I^c} A_i \left( x_i(\lbs) - [\xbs_{\LP}]^{(i)} \right) \nonumber \\
			&\leq b - \rhob + p \max_{i=1,\dots,m} \{ A_i x_i(\lbs) - A_i [\xbs_{\LP}]^{(i)} \}, \label{eq:feasibility_inequality}
	\end{align}
	where $I^c = \{1,\dots,m\}\setminus I$, and $b-\rhob$ constitutes an upper bound for $\sum_{i=1}^m A_i [\xbs_{\LP}]^{(i)}$ given that $\xbs_{\LP}$ is feasible for $\PbLP$.
	
	According to \cite[pag. 117]{shor1985minimization}, the component $[\xs_{\LP}]^{(i)}$ of the (unique, under Assumption~\ref{ass:uniqueness}) solution $\xbs_{\LP}$ to $\PbLP$ is the limit point of the sequence $\{\tilde{x}_i(k)\}_{k\geq 1}$, defined as
	\begin{equation*}
		\tilde{x}_i(k) = \frac{\sum_{r=1}^{k-1} \alpha(r)x_i(r+1)}{\sum_{r=1}^{k-1} \alpha(r)}.
	\end{equation*}
	By linearity, for all $k\geq 0$, we have that
	\begin{align*}
		A_i \tilde{x}_i(k) &= \frac{\sum_{r=1}^{k-1} \alpha(r) A_i x_i(r+1)}{\sum_{r=1}^{k-1} \alpha(r)} \\
			&\geq \min_{r\leq k} A_i x_i(r) \\
			&= \slb_i(k) \\
			&\geq \slb_i,
	\end{align*}
	where the first inequality is due to the fact that all $\alpha(k)$ are positive and the second equality follows from step~\ref{step:lower_bounds_update} of Algorithm~\ref{algo:Alg1}. In the final inequality, $\slb_i(k)$ is lower bounded by $\slb_i$, that denotes the limiting value of the non-increasing finite-valued sequence $\{\slb_i(k)\}_{k\geq 0}$. Note that all inequalities have to be intended component-wise.
	By taking the limit for $k\to\infty$, we also have that
	\begin{equation} \label{eq:xtilde_limit_lower_bound}
		A_i [\xbs_{\LP}]^{(i)} \geq \slb_i.
	\end{equation}
	
	By Proposition~\ref{prop:finite_time_primal_set}, there exists a finite iteration index $K$ such that $x_i(k)$ satisfies \eqref{eq:primal_solution_afterKeps}. Since \eqref{eq:feasibility_inequality} holds for any choice of $x_i(\lbs)$ which minimizes $(c_i\T + \lbs\!\T A_i)x_i$ over $\vert(X_i)$, if $k\geq K$, then we can choose $x_i(\lbs) = x_i(k)$. Therefore, for all $k\geq K$, \eqref{eq:feasibility_inequality} becomes
	\begin{align}
		\sum_{i=1}^m &A_i x_i(k) \nonumber \\
			&\leq b - \rhob + p \max_{i=1,\dots,m} \{ A_i x_i(k) - A_i [\xs_{\LP}]^{(i)} \} \nonumber \\
			&\leq b - \rhob + p \max_{i=1,\dots,m} \left\{ \max_{r\leq k} A_i x_i(r) - A_i [\xs_{\LP}]^{(i)} \right\} \nonumber \\
			&= b - \rhob + p \max_{i=1,\dots,m} \left\{ \sub_i(k) - A_i [\xs_{\LP}]^{(i)} \right\} \nonumber \\
			&\leq b - \rhob + p \max_{i=1,\dots,m} \left\{ \sub_i - \slb_i \right\} \nonumber \\
			&= b, \label{eq:primal_feasibility}
	\end{align}
	where the second inequality is obtained by taking the maximum up to $k$, the first equality is due to step~\ref{step:upper_bounds_update} of Algorithm~\ref{algo:Alg1}, the third inequality is due to the fact that $\sub_i$ is the limiting value of the non-decreasing finite-valued sequence $\{\sub_i(k)\}_{k\geq 1}$ together with \eqref{eq:xtilde_limit_lower_bound}, and the last equality comes from the definition of $\rho(k)=p \max \{\rho_1(k), \dots, \rho_m(k) \}$ where $\rho_i(k)=\sub_i(k)-\slb_i(k)$.
	
	From \eqref{eq:primal_feasibility} we have that, for any $k\geq K$, the iterates $x_i(k)$, $i=1,\dots,m$, generated by Algorithm~\ref{algo:Alg1} provide a feasible solution for \ref{eq:primal_program}, thus concluding the proof.
\end{proofof}

\begin{proofof}{Theorem~\ref{thm:primal_performance}}
	Denote as $\Jsp$, $\Jslpb$, and $\Jslp$ the optimal cost of \ref{eq:primal_program}, $\PbLP$, and \ref{eq:LP_primal_program}, respectively. From Assumption~\ref{ass:boundedness} it follows that $\Jsp$, $\Jslpb$, and $\Jslp$ are finite.
	
	Consider the quantity $\sum_{i=1}^m c_i\T x_i(k) - \Jsp$. \\
	As in the proof of Theorem~3.3 in \cite{vujanic2016decomposition}, we add and subtract $\Jslpb$ and $\Jslp$ to obtain
	\begin{align}
		\sum_{i=1}^m c_i\T x_i(k) - &\Jsp = \Big( \sum_{i=1}^m c_i\T x_i(k) - \Jslpb \Big) \nonumber \\
			&+ (\Jslpb - \Jslp) + (\Jslp-\Jsp).\label{eq:cost_exansion}
	\end{align}
	We shall next derive a bound for each term in \eqref{eq:cost_exansion}.

	\emph{Bound on $\sum_{i=1}^m c_i\T x_i(k) - \Jslpb$}:\\
	Similarly to the proof of Theorem~\ref{thm:primal_feasibility} for feasibility, due to Theorem~2.5 in \cite{vujanic2016decomposition}, have that there exists a set $I$ of cardinality at least $m-p$ such that $x_i(\lbs) = [\xbs_{\LP}]^{(i)}$, for all $i\in I$. Therefore,
	\begin{align}
		\sum_{i=1}^m c_i\T x_i&(\lbs) - \Jslpb \nonumber \\
			&= \sum_{i=1}^m c_i\T x_i(\lbs) - \sum_{i=1}^m c_i\T [\xbs_{\LP}]^{(i)} \nonumber \\
			&= \sum_{i\in I^c} c_i\T x_i(\lbs) - c_i\T [\xbs_{\LP}]^{(i)} \nonumber \\
			&\leq p \max_{i=1,\dots,m} \left\{ c_i\T x_i(\lbs) - c_i\T [\xbs_{\LP}]^{(i)} \right\}, \label{eq:cost_inequality}
	\end{align}
	where $I^c = \{1,\dots,m\}\setminus I$.
	
	According to \cite[pag. 117]{shor1985minimization}, the components $[\xbs_{\LP}]^{(i)}$ of the (unique, under Assumption~\ref{ass:uniqueness}) solution $\xbs_{\LP}$ to $\PbLP$ is the limit point of the sequence $\{\tilde{x}_i(k)\}_{k\geq 1}$, defined as
	\begin{equation*}
		\tilde{x}_i(k) = \frac{\sum_{r=1}^{k-1} \alpha(r)x_i(r+1)}{\sum_{r=1}^{k-1} \alpha(r)}.
	\end{equation*}
	By linearity, for all $k\geq 1$, we have that
	\begin{align*}
		c_i\T \tilde{x}_i(k) &= \frac{\sum_{r=1}^{k-1} \alpha(r) c_i\T x_i(r+1)}{\sum_{r=1}^{k-1} \alpha(r)} \geq \min_{r\leq k} c_i\T x_i(r) \geq \gammalb_i,
	\end{align*}
	where the first inequality is due to the fact that all $\alpha(k)$ are positive and the last one derives from the fact $\{\min_{r\leq k} c_i\T x_i(r)\}_{k\geq1}$ is a non-increasing sequence that takes values in a finite set, and hence is lower bounded by its limiting value $\gammalb_i$. Therefore, by taking the limit for $k\to\infty$, we also have that
	\begin{equation} \label{eq:xtilde_limit_cost_lower_bound}
		c_i\T [\xbs_{\LP}]^{(i)} \geq \gammalb_i.
	\end{equation}
	
	Since \eqref{eq:cost_inequality} holds for any choice of $x_i(\lbs)$ which minimize $(c_i\T + \lbs\!\T A_i)x_i$ over $\vert(X_i)$, by Proposition \ref{prop:finite_time_primal_set} it follows that, for $k\geq \bar K$, $x_i(\lbs) = x_i(k)$ and, as a result
	\begin{align}
		\sum_{i=1}^m c_i\T &x_i(k) - \Jslpb \nonumber \\
			&\leq p \max_{i=1,\dots,m} \left\{ c_i\T x_i(k) - c_i\T [\xbs_{\LP}]^{(i)} \right\} \nonumber \\
			&\leq p \max_{i=1,\dots,m} \left\{ \max_{r\leq k} c_i\T x_i(r) - c_i\T [\xbs_{\LP}]^{(i)} \right\} \nonumber \\
			&\leq p \max_{i=1,\dots,m} \left\{ \max_{r\leq k} c_i\T x_i(r) - \gammalb_i \right\}, \nonumber
	\end{align}
	where the second inequality is obtained by taking the maximum up to iteration $k$ and the third inequality is due to \eqref{eq:xtilde_limit_cost_lower_bound}.

	Now if we recall the definition of $\gamma(k)$ in \eqref{eq:gamma_update} and its finite-time convergence to $\gammaub$, jointly with the fact that $\gammalb_i$ is the limiting value of $\{\min_{r\leq k} c_i\T x_i(r)\}_{k\geq1}$, we finally get that there exists $K \ge \bar K$, such that for $k\ge K$
	\begin{equation*}
		p \max_{i=1,\dots,m} \left\{ \max_{r\leq k} c_i\T x_i(r) - \gammalb_i \right\}= \gammaub,
	\end{equation*}
	thus leading to
	\begin{align*}
		\sum_{i=1}^m c_i\T &x_i(k) - \Jslpb \leq \gammaub, \, k\ge K. 
	\end{align*}	

	\emph{Bound on $\Jslpb - \Jslp$:}\\
	Problem \ref{eq:LP_primal_program} can be considered as a perturbed version of $\PbLP$, since the coupling constraint of $\PbLP$ is given by
	\begin{equation*}
		\sum_{i=1}^m A_i x_i \leq b-\rhob
	\end{equation*}
	and that of \ref{eq:LP_primal_program} can be obtained by adding $\rhob$ to its right-hand-side.
	From perturbation theory (see \cite[Section~5.6.2]{boyd2004convex}) it then follows that the optimal cost $\Jslp$ is related to $\Jslpb$ by:
	\begin{equation} \label{eq:perturbation_theory}
		\Jslpb - \Jslp \leq \lbs\!\T \rhob.
	\end{equation}
	From Assumption~\ref{ass:slater}, by applying \cite[Lemma~1]{nedic2009approximate} we have that for all $\lambda\geq0$
	\begin{align}
		\norm{\lbs}_1 &\leq \frac{1}{m\zeta} \left( \sum_{i=1}^m c_i\T \hat{x}_i + \lambda\T b - \sum_{i=1}^m \min_{x_i\in X_i} (c_i\T + \lambda\T A_i)x_i \right) \nonumber \\
			&\leq \frac{1}{m\zeta} \left( \sum_{i=1}^m c_i\T \hat{x}_i - \sum_{i=1}^m \min_{x_i\in X_i} c_i\T x_i \right) \nonumber \\
			&\leq \frac{1}{\zeta} \max_{i=1,\dots,m} \left\{ \max_{x_i\in X_i} c_i\T x_i - \min_{x_i\in X_i} c_i\T x_i \right\} \nonumber \\
			&= \frac{\gammat}{p\zeta}, \label{eq:multipliers_boundedness}
	\end{align}
	where the second inequality is obtained setting $\lambda = 0$, the third inequality comes from the fact that $c_i\T \hat{x}_i \leq \max_{x_i\in X_i} c_i\T x_i$ and that $\sum_{i=1}^m \beta_i \leq m\max_i \beta_i$, and the third equality is due to \eqref{eq:paul_gamma}.
	Using \eqref{eq:multipliers_boundedness} in \eqref{eq:perturbation_theory} we have
	\begin{align*} \label{eq:perturbation_theory}
		\Jslpb - \Jslp &\leq \lbs\!\T \rhob \nonumber \\
			&\leq \norm{\lbs}_1 \norm{\rhob}_\infty \nonumber \\
			&\leq \frac{\norm{\rhob}_\infty}{p\zeta} \gammat, \nonumber
	\end{align*}
	where the second inequality is due to the H\"{o}lder's inequality.
	
	\emph{Bound on $\Jslp-\Jsp$:}\\
	Since \ref{eq:LP_primal_program} is a relaxed version of \ref{eq:primal_program}, then $\Jslp-\Jsp \leq 0$.
	
	The proof is concluded considering \eqref{eq:cost_exansion} and inserting the bounds obtained for the three terms.
\end{proofof}

\section{Application to optimal PEVs charging} \label{sec:exaple}

In this section we show the efficacy of the proposed approach in comparison to the one described in \cite{vujanic2016decomposition} on the Plug-in Electric Vehicles (PEVs) charging problem described in \cite{vujanic2016decomposition}. This problem consists in finding an optimal overnight charging schedule for a fleet of $m$ vehicles, which has to satisfy both local requirements and limitations (e.g., maximum charging power and desired final state of charge for each vehicle), and some network-wide constraints (i.e., maximum power that the network can deliver at each time slot). We consider both version of the PEVs charging problem, namely, the ``charge only'' setup in which all vehicles can only draw energy from the network, and the ``vehicle to grid'' setup where the vehicles are also allowed to inject energy in the network.

The improvement of our approach with respect to that in \cite{vujanic2016decomposition} is measured in terms of the following two relative indices: the reduction in the level of conservativeness
\begin{equation*}
	\Delta\rho_\% = \frac{\norm{\rhot}_\infty - \norm{\rhob}_\infty}{\norm{\rhot}_\infty} \cdot 100
\end{equation*}
and the improvement in performance achieved by the primal solution
\begin{equation*}
	\Delta J_\% = \frac{J_{\rhot}-J_{\rhob}}{J_{\rhot}} \cdot 100,
\end{equation*}
where $J_{\rhot} = \sum_{i=1}^m c_i\T x_i(\ls_{\rhot})$ and $J_{\rhob} = \sum_{i=1}^m c_i\T x_i(\lbs)$. A positive value for these indices indicates that our approach is less conservative.

For a thorough comparison we determined the two indices while varying: i) the number of vehicles in the network, ii) the realizations of the random parameters entering the system description (cost of the electrical energy and local constraints), and iii) the right hand side of the joint constraints. All parameters and their probability distributions were taken from \cite[Table~1]{vujanic2016decomposition}.

In Table~\ref{tab:results_agents} we report the conservativeness reduction and the cost improvement for the ``vehicle to grid'' setup. As it can be seen from the table, the level of conservativeness is reduced by $50\%$ while the improvement in performance (witnessed by positive values of $\Delta J_\% $) drops as the number of agents grows. This is due to the fact that the relative gap between $J_{\rhot}$ and $\Jsp$ tends to zero as $m\to\infty$, thus reducing the margin for performance improvement. 

\begin{table}[b]
	\centering
	\begin{tabular}{ccccccc}
		$m$				&	$250$		& $500$		& $1000$	& $2500$	& $5000$	& $10000$		\\ \hline
		$\Delta\rho_\%$	&	$50\%$		& $50\%$	& $50\%$	& $50\%$	& $50\%$	& $50\%$		\\ \hline
		$\Delta J_\%$	&	$13.9\%$	& $3.1\%$	& $1.1\%$	& $0.15\%$	& $0.05\%$	& $0.02\%$		\\ \hline
	\end{tabular}
	\caption{Reduction in the level of conservativeness ($\Delta\rho_\%$) and improvement in performance ($\Delta J_\%$) achieved by the primal solution obtained by the proposed method when compared with the one proposed in \cite{vujanic2016decomposition}.}
	\label{tab:results_agents}
\end{table}
We do not report the results for the ``charge only'' setup since the two methods lead to the same level of conservativeness and performance of the primal solution.

\begin{figure}[t]
	\centering
	\includegraphics[width=\columnwidth]{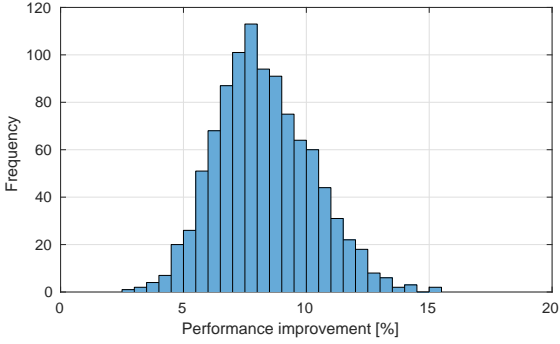}
	\caption{Histogram of the performance improvement ($\Delta J_\%$) achieved by the primal solution obtained by the proposed method with respect to the one proposed in \cite{vujanic2016decomposition} over $1000$ runs.}
	\label{fig:results_random}
\end{figure}
We also tested the proposed approach against changes of the random parameters defining the problem. We fixed $m = 250$ and performed $1000$ tests running Algorithm~\ref{algo:Alg1} and the approach in \cite{vujanic2016decomposition} with different realization for all parameters, extracted independently.
Figure~\ref{fig:results_random} plots an histogram of the values obtained for $\Delta J_\%$ in the $1000$ tests. Note that the cost improvement ranges from $3\%$ to $15\%$ and, accordingly to the theory, is always non-negative. The reduction in the level of conservativeness is also in this case $50\%$, suggesting that the proposed iterative scheme exploits some structure in the PEVs charging problem that the approach in \cite{vujanic2016decomposition} oversees.
Also in this case, in the ``charge only'' setup the two methods lead to the same level of conservativeness and performance.

Finally, we compared the two approaches in the ``vehicle to grid'' setup against changes in the joint constraints. If the number of electric vehicles is $m = 250$ and we decrease the maximum power that the network can deliver by $37\%$, then the $\rhot$ that results from applying the approach in \cite{vujanic2016decomposition} makes \ref{eq:LP_primal_program_rho} with $\rho = \rhot$ infeasible, thus violating Assumption~\ref{ass:uniqueness_paul}. Whereas with our approach \ref{eq:LP_primal_program_rho} with $\rho = \rhob$ remains feasible, $\rhob$ being the limiting value for $\{\rho(k)\}_{k\geq1}$ in Algorithm~\ref{algo:Alg1}.

\subsection{Performance-oriented variant of Algorithm~\ref{algo:Alg1}}
While Algorithm~\ref{algo:Alg1} is able to find a feasible solution to \ref{eq:primal_program}, it does not directly consider the performance of the solution, whereas the user is concerned with both feasibility and performance with higher priority given to feasibility. This calls for a modification to Algorithm~\ref{algo:Alg1} which also takes into account the performance achieved.

\begin{algorithm}[t]
	\begin{algorithmic}[1]
		\STATE {\it \% Initialize variables}
		\STATE $\lambda \gets 0$, $\sub_i \gets -\infty$, $\slb_i \gets +\infty$, $i=1,\dots,m$
		\STATE $\check{J} \gets +\infty$, $\delta \gets 0$, $k \gets 0$
		\REPEAT
			\FOR{$i=1$ \TO $m$}
				\STATE {\it \% Store tentative local solution}
				\IF{$\delta$ = 1}
					\STATE $\check{x}_i \gets x_i$
				\ENDIF
				\STATE {\it \% Update tentative local solution}
				\STATE $x_i \gets \arg\min_{x_i\in \vert(X_i)} (c_i\T + \lambda\T A_i)x_i$ 
			\ENDFOR
			\STATE {\it \% If solution is feasible and has better cost, then\\ ~~~~tell agents to update their tentative solutions}
			\IF{$\sum_{i=1}^m A_i x_i \leq b$ \AND $\sum_{i=1}^m c_i\T x_i < \check{J}$}
				\STATE $\check{J} \gets \sum_{i=1}^m c_i\T x_i$
				\STATE $\delta \gets 1$
			\ELSE
				\STATE $\delta \gets 0$
			\ENDIF
			\STATE {\it \% Update tightening}
			\STATE $\sub_i \gets \max \{ \sub_i, A_i x_i \}$, $i=1,\dots,m$ 
			\STATE $\slb_i \gets \min \{ \slb_i, A_i x_i \}$, $i=1,\dots,m$ 
			\STATE $\deltas_i \gets \sub_i - \slb_i$, $i=1,\dots,m$ 
			\STATE $\rho \gets \displaystyle p \max \{ \deltas_1, \dots, \deltas_m \}$ 
			\STATE {\it \% Update dual variables}
			\STATE $\lambda \gets [\lambda + \alpha(k) (\sum_{i=1}^m A_i x_i - b + \rho) ]_+$ 
			\STATE {\it \% Update iteration counter}
			\STATE $k \gets k+1$
		\UNTIL{time is over}
	\end{algorithmic}
	\caption{Performance-oriented version}
	\label{algo:Alg2}
\end{algorithm}

Theorem~\ref{thm:primal_feasibility} guarantees that there exists an iteration index $K$ after which the iterates stay feasible for \ref{eq:primal_program} for all $k\geq K$. Now, suppose that the agents, together with the $A_i x_i(k)$ also transmit $c_i\T x_i(k)$ to the central unit, then the central unit can construct the cost of $x(k) = [x_1(k)\T,\cdots,x_m(k)\T]\T$ at each iteration. When a feasible solution is found, its cost may be compared with that of a previously stored solution, and the central unit can decide to keep the new tentative solution or discard it. This way we are able to track the best feasible solution across iterations.

The modified procedure is summarized in Algorithm~\ref{algo:Alg2}. Note that, compared to Algorithm~\ref{algo:Alg1}, each agent is required to transmit also the cost of its tentative solution.

To show the benefits of Algorithm~\ref{algo:Alg2} in terms of performance, we run $1000$ test with $m = 250$ vehicles in the ``charge only'' setup, where we are also able to compute the optimal solution of \ref{eq:primal_program}, and compare the performance of Algorithm~\ref{algo:Alg1} and~\ref{algo:Alg2} in terms of relative distance from the optimal cost $\Jsp$ of \ref{eq:primal_program}.

\begin{figure}[t]
	\centering
	\includegraphics[width=\columnwidth]{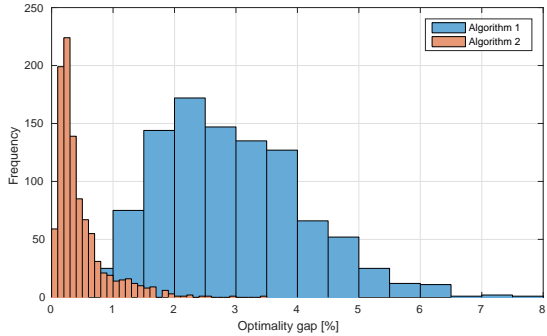}
	\caption{Histogram of the relative distance from the optimal value of \ref{eq:primal_program} achieved by the primal solution obtained by Algorithm~\ref{algo:Alg1} (blue) and Algorithm~\ref{algo:Alg2} (orange), over $1000$ runs.}
	\label{fig:results_random_alg2}
\end{figure}

Figure~\ref{fig:results_random_alg2} shows the distribution of $(J_{\rhob}-\Jsp)/\Jsp\cdot100$ obtained with Algorithm~\ref{algo:Alg1} (blue) and $(\check{J}-\Jsp)/\Jsp\cdot100$ obtained with Algorithm~\ref{algo:Alg2} (orange) for the $1000$ runs. As can be seen from the picture, most runs of Algorithm~\ref{algo:Alg2} result in a performance very close to the optimal one, while the runs from Algorithm~\ref{algo:Alg1} exhibit lower performance.

\section{Concluding remarks}
We proposed a new method for computing a feasible solution to a large-scale mixed integer linear program via a decentralized iterative scheme that decomposes the program in smaller ones and has the additional beneficial side-effect of preserving privacy of the local information if the problem originates from a multi-agent system.

This work improves over existing state-of-the-art results in that feasibility is achieved in a finite number of iterations and the decentralized solution is accompanied by a less conservative performance certificate. The application to a plug-in electric vehicles optimal charging problem verifies the improvement gained in terms of performance.

Future research directions include the development a distributed algorithm, which does not require any central authority but only communications between neighboring agents, and allows for time-varying communications among agents.

Moreover, we aim at exploiting the analysis of \cite{udell2016bounding} to generalize our results to problems with nonconvex objective functions.

\bibliographystyle{abbrv}

\end{document}